\pdfoutput=1
\documentclass[letterpaper, 10 pt, conference]{ieeeconf}  % Comment this line out
\IEEEoverridecommandlockouts                              % This command is only
\usepackage{graphics} % for pdf, bitmapped graphics files
\usepackage{epsfig} % for postscript graphics files
\usepackage{times} % assumes new font selection scheme installed
\usepackage{amsmath} % assumes amsmath package installed
\usepackage{amssymb}  % assumes amsmath package installed
\usepackage{bbm}
\usepackage{color}
\newtheorem{thm}{Theorem}

\title{\LARGE \bf
Optimal Power Flow with Weighted Chance Constraints \\ and General Policies for Generation Control
}

\author{Line Roald, Sidhant Misra, Michael Chertkov, G\"oran Andersson% <-this % stops a space
\thanks{L. Roald and G. Andersson are with the Power Systems Laboratory, Department of Electrical Engineering,
        ETH Zurich, Switzerland
        {\tt\small \{roald|andersson\}@eeh.ee.ethz.ch}}%
\thanks{S. Misra and M. Chertkov. are with the Center for Nonlinear Studies and Theoretical Division T-4 of Los Alamos National Laboratory,
        Los Alamos, NM 87544. M. Chertkov is also affiliated with the New Mexico Consortium, Los Alamos, NM 87544.
        {\tt\small \{sidhant|chertkov\}@lanl.gov}}
}

\begin{document}

\maketitle
\thispagestyle{empty}
\pagestyle{empty}

%%%%%%%%%%%%%%%%%%%%%%%%%%%%%%%%%%%%%%%%%%%%%%%%%%%%%%%%%%%%%%%%%%%%%%%%%%%%%%%%
\begin{abstract}
Due to the increasing amount of electricity generated from renewable sources, uncertainty in power system operation will grow. This has implications for tools such as Optimal Power Flow (OPF), an optimization problem widely used in power system operations and planning, which should be adjusted to account for this uncertainty. One way to handle the uncertainty is to formulate a Chance Constrained OPF (CC-OPF) which limits the probability of constraint violation to a predefined value.
However, existing CC-OPF formulations and solutions are not immune to drawbacks. On one hand, they only consider affine policies for generation control, which are not always realistic and may be sub-optimal. On the other hand, the standard CC-OPF formulations do not distinguish between large and small violations, although those might carry significantly different risk. In this paper, we introduce the Weighted CC-OPF (WCC-OPF) that can handle general control policies while preserving convexity and allowing for efficient computation. The weighted chance constraints account for the size of violations through a weighting function, which assigns a higher risk to a higher overloads. We prove that the problem remains convex for any convex weighting function, and for very general generation control policies. In a case study, we compare the performance of the new WCC-OPF and the standard CC-OPF and demonstrate that WCC-OPF effectively reduces the number of severe overloads. Furthermore, we compare an affine generation control policy with a more general policy, and show that the additional flexibility allow for a lower cost while maintaining the same level of risk.
\end{abstract}

%%%%%%%%%%%%%%%%%%%%%%%%%%%%%%%%%%%%%%%%%%%%%%%%%%%%%%%%%%%%%%%%%%%%%%%%%%%%%%%%
\section{Introduction}
Optimal Power Flow (OPF) is an optimization problem widely used in power system operational and expansion planning. In this paper we will address the application of OPF in operational planning. The objective is to minimize cost of operation, while avoiding violations of technical constraints such as limits on transmission or generation capacity.
However, with increasing penetration of electricity generation from renewable sources, power systems operators face higher operational uncertainty. To avoid violations of technical constraints, it is increasingly important to account for these uncertainties during the operational planning phase, i.e., within the OPF problem.

One way to handle uncertainty is to limit the probability of constraint violation by formulating a chance constrained OPF (CC-OPF).
In \cite{maria}, a joint chance constrained problem was formulated, and solved using the scenario approach, a solution method based on samples.
Both \cite{12BCH}, \cite{line} formulated the CC-OPF with separate chance constraints, and obtained an analytical reformulation assuming normally distributed wind power fluctuations.
Convexity of the CC-OPF with respect to the characteristics of controllable generation was a key observation of \cite{12BCH} which also allowed for guarantees of exactness and lead to efficient implementation for large networks (thousands of nodes).

While previous versions of CC-OPF have shown the ability to limit the probability of violations, they have some drawbacks.
First, current CC-OPF formulations only consider affine policies for generation control, which is not necessarily realistic for large wind power deviations and may lead to sub-optimal results.
Second, standard chance constraints do not distinguish between large and small violations, although those might carry significantly different risk.
This paper addresses both drawbacks, thus extending existing formulations in two ways.

First, we suggest more general policies for generator control. An affine policy is a good approximation of the existing implementations of automatic generation control (AGC) for small deviations. However, a more general policy allows us to model additional control actions, such as manual activation of reserves during large wind power deviations. It is thus possible to obtain a cheaper, more realistic solution while maintaining the same level of risk. To the best of our knowledge, no previous chance constrained OPF aimed at reducing the cost in a similar way.

Second, the proposed extension accounts for both the probability and the size of constraint violations.
A standard chance constraint only considers the probability of violations, and does not account for the risk related to the actual constraint violations.
To remedy this, we define a weighting function for the overload, and calculate the risk as the product between the weighting function and the probability distribution of the overload.
This Weighted CC-OPF (WCC-OPF) is convex with respect to parameters of the generation, thus allowing exactness and tractability in the spirit of \cite{12BCH}.
While the set-up allows for general convex weighting functions, the considered examples are chosen to be applicable to the power system problem. The examples are closely related to the severity functions applied in risk-based OPF (e.g., \cite{Xiao2009}, \cite{Gabi2012}, \cite{ownPSCC}), where they are used to account for the severity of post-contingency overloads.
Such a severity function was combined with a CC-OPF in \cite{ownPSCC}. However, \cite{ownPSCC} limited the probability of high risk situations, and did not consider the average risk due to constraint violation.
The weighted chance constraint suggested here is somewhat similar to a constraint on the conditional-value-at-risk (CVaR), a risk measure borrowed from finance, which was introduced by \cite{rockafellar} and applied to the OPF problem in \cite{Tyler}. However, the weighted chance constraint suggested in this paper is more general and adapted to better fit the power system problem.

The remainder of the paper is organized as follows. Section II introduces the original CC-OPF. In Section III, the weighted chance constraint is introduced, along with a proof of convexity. Section IV provides some relevant examples for practical weighting functions with affine and piecewise affine generation control. Section V illustrates the performance of the method based on the IEEE RTS96 test system.
Section VI summarizes and concludes the paper.

\subsection{Notations}
We denote vectors by lower case letters $p, \omega$.
The components of the vectors are denoted by using subscripts, i.e, the $i$th component of $p$ is denoted by $p_i$.
Matrices are denoted by upper/lower bold case letters, $\boldsymbol{\alpha}, \mathbf{M}$, and $\boldsymbol{\alpha}_{(i,\cdot)}, \boldsymbol{\alpha}_{(\cdot,i)}$ denote the $i^{th}$ row and column of $\boldsymbol{\alpha}$, respectively.
Index $i$ refers to generators, while index $ij$ refers to lines.

\section{Chance constrained Optimal Power Flow}
In this section, we review the existing literature regarding the modeling of the OPF problem with fluctuating in-feeds as a CC-OPF. Our description is based on the formulation in \cite{12BCH}.
Let $\mathbbm{G} = (\mathcal{V}, \mathcal{E})$ represent the graph of the power network, where $\mathcal{V}$ is the set of nodes with $|\mathcal{V}| = m$ and $\mathcal{E}$ is the set of edges/lines of the system
with $|\mathcal{E}| = n$.
The set of wind generators is denoted by $\mathcal{W}\subseteq\mathcal{V}$. Wind is considered as the only source of fluctuations in this paper, although the formulation can easily be extended to handle fluctuations from other sources such as solar PV or load.
The set of non-wind generators is denoted by $\mathcal{G} \subseteq \mathcal{V}$, and are assumed to be controllable within their limits.
To simplify notation, we assume that there is one controllable generator $p$, one wind generator $w$ and one demand $d$ per node, such that $|\mathcal{G}|=|\mathcal{W}|=|\mathcal{V}| = m$. Nodes without generation or load can be handled by setting the respective entries to zero.

\subsubsection{Wind power in-feeds}
The wind in-feeds $w$ are modeled as the sum of the forecasted electricity production from wind, given by $v=\mathbb{E}_{w}[w]$, and a zero mean fluctuating component ${\omega}\doteq w-v$.
\subsubsection{Generation control}
Since secure operation of the power system requires balance between produced and consumed power at all times, any deviation $\omega$ in the wind power production must be balanced by an adjustment in the controllable generation. In the CC-OPF \cite{12BCH}, these adjustments are modelled through an affine policy, reflecting the automatic generation control which is establishing balance within tens of second to a few minutes \cite{wollenberg}:
\begin{equation}
\tilde{p}(\omega) = p-\boldsymbol{\alpha}\omega~,
\label{affine}
\end{equation}
where $p$ is the vector of scheduled generation, and $\boldsymbol{\alpha} \in \mathbbm{R}^{m \times m}$ is a matrix with elements describing the response of a generator to a wind fluctuation. The elements of each column in $\boldsymbol{\alpha}$ sum to one.%, $\sum\limits_{i\in\mathcal{G}} \boldsymbol{\alpha}_{(\cdot,i)}  =  1$.
If we assume that the generators respond only to the total deviation in wind generation output, as in \cite{12BCH}, all columns of $\boldsymbol{\alpha}$ are identical, i.e., $\boldsymbol{\alpha}_{(\cdot,i)} = \boldsymbol{\alpha}_{(\cdot,j)}~\forall_{\{i,j\}}$.  %Note that $\boldsymbol{\alpha}_{(i,\cdot)}$ is the $i^{th}$ column and $\boldsymbol{\alpha}_{(\cdot,i)}$ is the $i^{th}$ row of $\boldsymbol{\alpha}$.
\subsubsection{Power Flows}
The power flows $p_{ij}$ on each line are computed according to the standard so-called DC approximation \cite{wollenberg}
\begin{equation}
p_{ij}=\mathbf{M}_{(ij,\cdot)}(p-\boldsymbol{\alpha}\omega+v+{\omega}-d), \forall_{ij\in \mathcal{E}}.
\label{eq:lineflows}
\end{equation}
The matrix $\mathbf{M} \in \mathbbm{R}^{n \times m}$ relates the line flows to the nodal power injections, which are expressed as the sum of generation $p-\boldsymbol{\alpha}\omega$, wind power production $v+\omega$ and demand $-d$.
$M$ is defined as
\begin{equation}
\mathbf{M} = B_f \begin{bmatrix} (\widetilde{B}_{bus})^{-1}   ~~~ \bold{0}  \\ ~~\bold{0}  ~~~~~~~~~ 0\end{bmatrix}
\end{equation}
where ${B}_{f}$ is the line susceptance matrix and $\widetilde{B}_{bus}$ the bus susceptance matrix (without the last column and row) \cite{maria}. $\mathbf{M}_{(ij,\cdot)}$ is the row of $M$ related to the line $(ij) \in \mathcal{E}$.
\subsubsection{Optimization Problem}
With the above modelling considerations, the CC-OPF is stated as follows:
\begin{align}%{2}
& \min\limits_{p,\boldsymbol{\alpha}} \ \sum_{i\in\mathcal{G}} c_i p_i \label{CC-OPF} \\
& \textrm{s.t.} \\
& \sum_{i\in\cal{V}}p_i - d_i + v_i = 0 \label{powerbalance} \\
& \sum\limits_{i\in\mathcal{G}} \boldsymbol{\alpha}_{(i,j)}  =  1~~\forall_{j\in\mathcal{W}},   \quad p \ge 0, \quad  \boldsymbol{\alpha} \ge 0    \label{conic-first} \\
& \mathbb{P}_{\omega}\left[p_i - \boldsymbol{\alpha}_{(i,\cdot)}\omega > p^{max}_i \right] < \epsilon_i~~ \forall_{i \in \mathcal{G}} \label{chance_gen_min} \\
& \mathbb{P}_{\omega}\left[p_i - \boldsymbol{\alpha}_{(i,\cdot)}\omega < p^{min}_i \right] < \epsilon_i~~ \forall_{i \in \mathcal{G}} \label{chance_gen_max} \\
& \mathbb{P}_{\omega}\left[\mathbf{M}_{(ij,\cdot)}(p-\boldsymbol{\alpha}\omega-d+v+{\omega}) > ~\bar{p}_{ij}\right] <  \epsilon_{ij}~~~ \forall_{\{ij\} \in \mathcal{E}} \label{upperlinechance} \\
& \mathbb{P}_{\omega}\left[\mathbf{M}_{(ij,\cdot)}(p-\boldsymbol{\alpha}\omega-d+v+{\omega}) < -\bar{p}_{ij}\right] <  \epsilon_{ij}~~ \forall_{\{ij\} \in \mathcal{E}} \label{lowerlinechance}
\end{align}
The objective \eqref{CC-OPF} is to minimize generation cost, where $c$ contains costs, i.e. bids, from the generators.
Eq. \eqref{powerbalance} reflects the power balance during nominal (pre-planned) operation, while \eqref{conic-first} accounts for deviations from nominal, where $p$ and $\alpha$ are assumed positive. (The latter constraints is not theoretically significant, but reflects engineering practice).
Since the generation output and the power flows depend on the fluctuations $\omega$, the generation constraints \eqref{chance_gen_min}, \eqref{chance_gen_max} and transmission constraints \eqref{upperlinechance}, \eqref{lowerlinechance} are formulated as chance constraints with accepted violation probability $\epsilon$.
To obtain a tractable optimization problem, the chance constraints must be reformulated. Assuming that the wind fluctuations follow a multivariate normal distribution, \eqref{chance_gen_min} - \eqref{lowerlinechance} can be restated as second-order cone constraints \cite{12BCH}.

The CC-OPF formulation above accounts for uncertainty, but have some shortcomings. First, it assumes that the generation control follows an affine policy.
Second, the chance constraints treat all violations equally, irrespective of the size of the violation.
In the following, we discuss a new way of formulating the chance constraints to address both those shortcomings.

\section{Weighted Chance Constraints}
% Put general constraint, to explain what it looks like
To formulate a chance constraint that can handle general policies for generation control and at the same time differentiates between large and small violations, we propose a generalized weighted chance constraint of the form
\begin{equation}
\int_{-\infty}^{\infty}f(y(\omega))P(\omega)d\omega \leq \epsilon~.
\label{cvar_general_1}
\end{equation}
Here, $P(\omega)$ is the multivariate distribution function of the fluctuations. The quantity $y(\omega)$ denotes the magnitude of overload and depends on whether we are considering a violation of the upper or lower limit on generation or power flow. In case of violation of the upper generation limit for generator $i$, we have
\begin{align}
	y(\omega) = \tilde{p}_i(\omega) - p_{i}^{max},  \label{y_gen_u}
\end{align}
and in the case of lower limit violation, we have
\begin{align}
	y(\omega) = p_{i}^{min} - \tilde{p}_i(\omega).  \label{y_gen_l}
\end{align}
Similarly, for power flows we have that for line $(i,j) \in \mathcal{E}$,
\begin{align}
	y(\omega) = \tilde{p}_{ij}(\omega) - p_{ij}^{max},  \label{y_line_u}
\end{align}
and
\begin{align}
	y(\omega) = p_{ij}^{min} - \tilde{p}_{ij}(\omega),  \label{y_line_l}
\end{align}
where
\begin{align}
	\tilde{p}_{ij}(\omega) = (\mathbf{M}_{(ij,\cdot)}(\tilde{p}(\omega)-d+v+{\omega})). \label{line_power_general}
\end{align}

 %and $y(\omega) = \tilde{p}(\omega) - p_{max}$ denotes the overload of a generator or transmission line as a function of the fluctuations.
Whenever we have $y>0$, it indicates a violation of the limit, while $y<0$ implies that we are in a safe operating region. The weighting function $f(y(\omega))$, which is nonzero only if $y>0$, describes the risk related to the overload.
Note that the constraint (\ref{cvar_general_1}) depends implicitly on the generation control policy $\tilde{p}(\omega)$.

The weighting function $f(y)$ can be implemented in many ways. For example, \eqref{cvar_general_1} is equivalent to a standard chance constraint if $f(y)$ is the unit step function, i.e., 0 for $y<0$ and 1 for $y\geq 0$. %However, the step function, and thus the original chance constraint, has some drawbacks. First, it does not distinguish between constraint violations of different size.
%Particularly for transmission lines, where the amount of overload has a significant influence on the risk of failure, the indicator function might not be the most appropriate way of expressing the risk.
However, the step function is not convex, which means that the standard chance constraint will not always be a convex constraint. However, as we will see in Theorem \ref{thm:convexity}, the constraint
(\ref{cvar_general_1}) is a convex constraint for general generation control policies and general probability distributions of the wind fluctuations whenever the weight function $f(.)$ is convex.

\begin{figure}[!t]
%    \psfrag{116}[cc][][1][0]{116}
\includegraphics[width=0.45\textwidth]{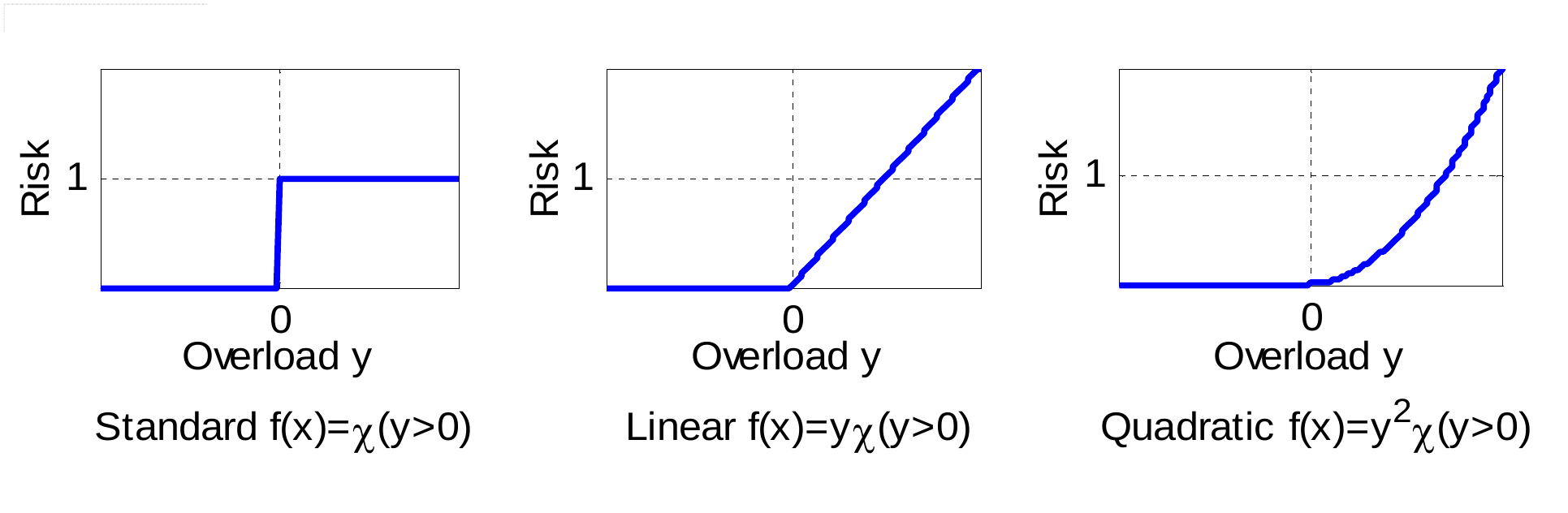}
\centering
\caption{Examples of weighting functions for chance constraints: Standard (left), linear (middle) and quadratic (right).}
\label{example1}
\end{figure}

While the desire to maintain convexity is a motivation on its own, a convex weight function, which assigns a higher risk to a higher overload, is also natural from an engineering point of view.
For example, in the context of so-called risk-based OPF, linear \cite{Xiao2009} and quadratic \cite{Gabi2012} weighting functions have been used to model the risk of post-contingency overloads.
Assuming that the risk of overload would be similar in the case of overloads induced by wind power fluctuations, weighting functions similar to these used in \cite{Xiao2009}, \cite{Gabi2012} could be considered.
For comparison, these weight functions are plotted in Fig.~\ref{example1} along with the weight function corresponding to a standard chance constraint (which is a unit step function).
While unit step function assigns equal risk to all magnitudes of overloads,
the risk increases linearly or quadratically with increasing overload for the linear or quadratic weight functions.

During large wind power deviations, the power system operator might use additional controls, such as slower reserves (e.g., spinning or tertiary reserves), to balance the system and control power flows.
Since these controls have a significant influence on the ability to handle disturbances, they should be included in the OPF to avoid unrealistic or sub-optimal solutions.
A main advantage of (\ref{cvar_general_1}) with a convex weighting function $f(y)$ is therefore in the principal ability to handle more general control policies while still maintaining convexity, as shown below.

\subsection{Convexity of weighted chance constraints}
We state our result regarding the convexity of the weighted chance constraint defined in (\ref{cvar_general_1}) in Theorem \ref{thm:convexity}. \\

\begin{thm} \label{thm:convexity}
	Let $\tilde{p}(\omega)$ be any general generator control policy.
	Consider a general weighted chance constraint defined in (\ref{cvar_general_1}) given by
	\begin{align}
		\int _{-\infty}^{\infty} f(y(\omega)) P(\omega) d\omega \leq \epsilon \label{cvar_general_repeat}
	\end{align}
	Assume that the weighting function $f(.)$ is a convex function of its argument. Then (\ref{cvar_general_repeat}) is a convex constraint in the control policy $\tilde{p}(\omega)$.
\end{thm}
\vspace{0.1in}
\begin{proof}
	Referring to Eqs.~(\ref{y_gen_u})-(\ref{y_line_l}) and Eq.~(\ref{line_power_general}), we see that the overload $y(\omega)$ is always a linear function of the generator control policy $\tilde{p}$.
	Combined with the convexity of the weight function $f(.)$, this implies that the function $f(y(\omega))$ is a convex function of $\tilde{p}$ for every $\omega$.
	The proof then follows from the fact that the expectation of a convex function is also convex.
\end{proof}

Theorem \ref{thm:convexity} shows that the weighted chance constraints are convex under general control policies, as long as the weighting function $f(.)$ is convex.
In addition, since the the network must obey total power balance at all times, we must also have the additional constraint on the generation policy that enforces this, given by
\begin{align}
	\sum_{i} \tilde{p}_i(\omega)  + v_i + \omega_i - d_i = 0, \quad \forall \omega. \label{total_power_balance}
\end{align}
Note that for every $\omega$, the constraint in (\ref{total_power_balance}) is simply a linear constraint.
This means that the control policy $\tilde{p}$ can be a variable within a convex WCC-OPF problem.
However for practical implementation of the weighted chance constraints within a standard convex optimization program, it is necessary to represent the control policy with finitely many parameters. It is straight forward to see that the convexity of (\ref{cvar_general_1}) and the linearity of (\ref{total_power_balance}) is preserved as long as we use any linear-in-parameter representation of the control policy $\tilde{p}$ given by
\begin{align}
	\tilde{p}(\omega) = p - \sum_{k = 1}^{K} \alpha_k g_k(\omega), \label{linear_in_param}
\end{align}
along with the constraints
\begin{align}
	 &p  + v - d = 0, \\
	 \sum_{k = 1}^{K} &\alpha_k g_k(\omega)  + \omega = 0,
\end{align}
where $\alpha_k \in \mathbbm{R}^m$ and $g_k(\omega)\in \mathbbm{R}$.
The original affine policy \eqref{affine} can be obtained as a special case of (\ref{linear_in_param}) by setting $K = m$ and $g_k(\omega) = \omega_k$.

In the above, $g_k(\omega)$ can be any general function of $\omega$. The expression for the control policy in (\ref{linear_in_param}) has extensive representation potential. For example, we can represent most general
control policies $\tilde{p}(\omega)$ (e.g., continuous or piece-wise continuous) approximately in the form of (\ref{linear_in_param}) by using a piece-wise constant representation of $\tilde{p}(\omega)$, i.e.,
\begin{align}
	\tilde{p}(\omega) = p - \sum_{k = 1}^{K} \alpha_k \chi(\omega \in S_k),
\end{align}
where $S_k$ are sets that form a partition of the domain of $\omega$ and $\alpha_k$ is the value of $\tilde{p}$ on $S_k$.
Increasing the number of terms $K$ in the representation allows us to approximate general control policies with higher fidelity,
but the number of variables to optimize over in the WCC-OPF also increases.
%but at the same time there is a trade-off because the number of variables to optimize over in the WCC-OPF also increases.
Thus, a trade-off between fidelity and computational effort must be made.

\section{Weighted chance constraints with example weight functions and control policies for normally distributed fluctuations}

%\section{Weighted chance constraints with example weighting functions and affine control}
In this section, we derive the expressions for the weighted chance constraints for linear and quadratic weight functions and for affine and piece-wise affine control policies. The expressions assume that the
wind fluctuation $\omega$ is distributed as a multivariate normal distribution. For completeness, we also contrast these expressions with the original chance constraint in our current notation.

\subsection{Linear weight function and affine policy}
For the linear weight function, substituting $f(y) = y \chi(y > 0)$ in (\ref{cvar_general_1}), we get
\begin{equation}
\int_{-\infty}^{\infty}y\chi(y>0)P(\omega)d\omega~= \int_{0}^{\infty}y~P(y)dy \leq \epsilon. \label{linear}
\end{equation}
%which is the expectation of the truncated normal distribution with lower bound 0.

We derive the expressions for constraint violations of the upper generator and line limits. The lower limit case can be handled similarly.
Substituting $\tilde{p}(\omega) =  p - \boldsymbol{\alpha} \omega$ in (\ref{y_gen_u}), we get that the generator overloads $y_{i}^{u}$ are distributed as Gaussian random variables
with the average overloads $\mu_{i}^u$ and variances $(\sigma_{i}^u)^2$ given by
\begin{equation}
\mu_{i}^u  =  p_i - p^{max}_i~, \quad (\sigma_{i}^u)^2 =  \boldsymbol{\alpha}_{(i,\cdot)} \Sigma \boldsymbol{\alpha}_{(i,\cdot)}^T~. \label{mean_gen_affine}
\end{equation}
Similarly, using (\ref{y_line_u}), we get that the line overload $y_{ij}^u$ is Gaussian with
\begin{align}
\mu_{ij}^u &= \mathbf{M}_{(ij,\cdot)}(p-d+\mu) - \bar{p}_{ij}~, \\
(\sigma_{ij}^u)^2 &= \mathbf{M}_{(ij,\cdot)}(E-\boldsymbol{\alpha}) \Sigma (M_{(ij,\cdot)}(E-\boldsymbol{\alpha}))^T~. \label{variance_line_affine}
\end{align}
Similar expressions can be obtained for the lower limits $y_{i}^{l},~y_{ij}^l$ using (\ref{y_gen_l}) and (\ref{y_line_l}).
We observe that $y$ is distributed as a Gaussian random variable in all cases, since they are linear combinations of the random Gaussian $\omega$.

Since $y$ is a Gaussian random variable, the LHS of (\ref{linear}) is simply the expectation of a truncated Gaussian random variable, and can be expressed as
\begin{equation}
\mu \left(1 - \Phi\left( \frac{-\mu}{\sigma} \right)\right) + \frac{\sigma}{\sqrt{2\pi}}e^{-\frac{1}{2}\left(\frac{-\mu}{\sigma}\right)^2} \leq \epsilon~, \label{lin}
\end{equation}
where $\Phi(x)$ is the cumulative distribution function of the standard Gaussian. The weighted chance constraints for generator and line overloads are obtained by using the corresponding means and variances from Eqs.~(\ref{mean_gen_affine})-(\ref{variance_line_affine}) above in (\ref{lin}).
%and (\ref{y_line_u}), we

\subsection{Quadratic weight function and affine policy}
Using
$f(y) = y^2 \chi(y > 0)$ in (\ref{cvar_general_1}), we get
\begin{equation}
\int_{-\infty}^{\infty}y^2\chi(y>0)P(\omega)d\omega~= \int_{0}^{\infty}y^2~P(y)dy \leq \epsilon. \label{quadratic}
\end{equation}
The LHS of (\ref{quadratic}) is the second moment of a truncated normal and can be rewritten as
\begin{equation}
(\mu^2 + \sigma^2) \left(1 - \Phi\left( \frac{-\mu}{\sigma} \right)\right) + \frac{\mu\sigma}{\sqrt{2\pi}}e^{-\frac{1}{2}\left(\frac{-\mu}{\sigma}\right)^2} \leq \epsilon~. \label{quad}
\end{equation}
As for the case of a linear weight function, the weighted chance constraints for the quadratic weight function can be obtained by substituting the means and variances from Eqs.~(\ref{mean_gen_affine})-(\ref{variance_line_affine})
in (\ref{quad}).

\subsection{Standard chance constraints and affine policy}
Here we include the expressions for the standard chance constraints with $f(y)=\chi(y>0)$ for the sake of completeness.

\begin{align}
\int_{-\infty}^{\infty}\chi(y>0)P(\omega)d\omega~=\int_{0}^{\infty}P(y)dy \leq \epsilon
\end{align}
which can be reformulated as
\begin{equation}
\mu + \Phi^{-1}(1-\epsilon)\sigma \leq 0, \label{cc}
\end{equation}
and then the final form can be obtained by using Eqs.~(\ref{mean_gen_affine})-(\ref{variance_line_affine})  as usual.

\subsection{Linear weight function and piece-wise affine policy} \label{sec:lin_piece}
In this section, we consider a piece-wise affine policy which %, which is a more realistic policy taking into account activation of spinning reserves during large wind deviations.
%This policy
can be expressed as
\begin{equation}
\tilde{p}\doteq
\begin{cases}
p - \boldsymbol{\alpha}\omega,\quad & \Omega^-\leq \Omega\leq {\Omega}^+\\
p - \boldsymbol{\alpha}\omega + \beta^+_i,\quad & \Omega\geq \Omega^+>0\\
p - \boldsymbol{\alpha}\omega + \beta^-_i,\quad & \Omega\leq \Omega^-<0
\end{cases}
\label{depPol}
\end{equation}
Here, $\Omega \doteq \sum_{i\in \mathcal{V}} \omega_i$ is the total wind power deviation, and $\beta^+,~\beta^-$ represent additionally deployed reserves in case of large wind fluctuations, when $\Omega$ is larger (or smaller) than a given threshold $\Omega^+$ (or $\Omega^-$). To ensure a balanced system, we additionally enforce $\sum_{i\in\mathcal{G}} \beta^+ = \sum_{i\in\mathcal{G}} \beta^- = 0$.

Under a piece-wise affine policy, the linearly weighted chance constraint can be written as
\begin{align}
	&\int_{0}^{\infty} y P (y) dy = \int_{-\infty}^{\infty} \int_{0}^{\infty} y P(y|\Omega) P(\Omega) dy d\Omega  \label{linear_piece1} \\
	=& \int_{-\infty}^{\Omega^-} \int_{0}^{\infty} y P(y|\Omega) P(\Omega) dy d\Omega \nonumber \\
    +& \int_{\Omega^-}^{\Omega^+} \int_{0}^{\infty} y P(y|\Omega) P(\Omega) dy d\Omega \nonumber \\% \label{linear_piece2}
	+& \int_{\Omega^+}^{\infty} \int_{0}^{\infty} y P(y|\Omega) P(\Omega) dy d\Omega. \label{linear_piece3}
\end{align}
The random variable $y|\Omega$ is normally distributed and its mean and variance depends on which of the three regions in (\ref{depPol}) the total wind deviation $\Omega$ belongs to.
For generator power, the conditional mean and variance are given by
\begin{align}
	\mu_i^u(\Omega)  = p_i - p_i^{max} + \beta_i(\Omega) - \left(\mathbf{1}^T \Sigma \mathbf{1} \right)^{-1}  \left(\boldsymbol{\alpha}_{(i,.)} \Sigma \mathbf{1} \right) \Omega,
\label{piecewiseaffine}
\end{align}
where
\begin{align}
	\beta_i(\Omega) = \begin{cases} \label{beta_cases}
			\beta_i^-, \quad & \Omega < \Omega^-, \\
			0, \quad & \Omega^- \leq \Omega \leq \Omega^+, \\
			\beta_i^+, \quad & \Omega > \Omega^+.
		\end{cases}
\end{align}
and
\begin{align}
	(\sigma_i^{u})^2 = \boldsymbol{\alpha}_{(i,.)} \Sigma \boldsymbol{\alpha}_{(i,.)}^T - \left(\mathbf{1}^T \Sigma \mathbf{1} \right)^{-1}  \left(\boldsymbol{\alpha}_{(i,.)} \Sigma \mathbf{1} \right)^2.
\end{align}

Similarly, we can calculate the conditional means and variances of the line overloads as
\begin{align}
	\mu_{ij}^{u}(\Omega) &= \mathbf{M}_{(ij,.)} \left[  p - d + v + \beta(\Omega) \right] - p_{ij}^{max} \nonumber \\
                         &+ \left(\mathbf{1}^T \Sigma \mathbf{1} \right)^{-1} \left(\mathbf{M}_{(ij,.)} (I-\alpha) \Sigma \mathbf{1}\right) \Omega\label{moments_piece1}  \\
	(\sigma_{ij}^u)^2 &= \mathbf{M}_{(ij,.)}(I - \boldsymbol{\alpha})\Sigma(I-\boldsymbol{\alpha})^T \mathbf{M}_{(ij,.)}^T \nonumber\\
	& - \left(\mathbf{1}^T \Sigma \mathbf{1} \right)^{-1} \left( \mathbf{M}_{(ij,.)}(I - \boldsymbol{\alpha}) \Sigma \mathbf{1}  \right)^2, \label{moments_piece3}
\end{align}
with $\beta(\Omega)$ defined as in (\ref{beta_cases}).

The final version of the linear weighted chance constraints is obtained by substituting the expression for the expectation of a truncated Gaussian from (\ref{lin}) in (\ref{linear_piece3}), and insert the means and variances from (\ref{moments_piece1}) - (\ref{moments_piece3}).

\subsection{Quadratic weight function and piece-wise affine policy}
Following the same derivation as in Section \ref{sec:lin_piece}, we can obtain an expression analogous to (\ref{linear_piece3}) by simply replacing $y$ by $y^2$. The rest of the
derivation is identical with the only difference that we substitute the expression for the second moment of a truncated Gaussian from (\ref{quad}).

\section{Case Study}
The WCC-OPF is illustrated in a case study based on the IEEE RTS96 system, with modifications to the base case \cite{RTS96} similar to these described in \cite{line}. The line capacities are reduced to 80 \% of the steady state capacity listed in \cite{RTS96}, since only active power flows are considered. The generation costs used in the OPF are linear costs adopted from \cite{milano}, and we also assume that the generation capacity is twice larger than those from the original RTS96 model \cite{RTS96}. The minimal generator output is set to 0 for all generators. Aggregations of wind power plants are located at bus 8 and bus 15, with a forecasted in-feed of 125 and 175 MW, respectively. The wind power fluctuations are assumed to follow a multivariate normal distribution with zero mean and standard deviations of 9.4 and 13.1 MW, respectively, with a correlation coefficient $\rho=0.2$. Finally, we assume that the generation policy only reacts to the overall wind fluctuation $\Omega$, such that \eqref{affine} reduces to
\begin{equation}
\tilde{p}(\omega) = p-\alpha\Omega~.
%\label{affine}
\end{equation}
Here, $\alpha \in \mathbb{R}^m$ is a vector with one contribution factor for each generator.

In the following, we investigate two cases. First, we assume an affine policy for generation control, and assess how the weighted chance constraints perform compared with standard chance constraints in terms of number and size of the constraint violations. Second, we investigate how a the more flexible piece-wise affine policy for generation control \eqref{depPol} performs compared with the standard affine policy \eqref{affine} in the case of the linear WCC-OPF.

\subsection{Comparison of CC-OPF and WCC-OPF}
We compare the standard CC-OPF with the linear and quadratic WCC-OPF assuming an affine policy for generation control.
We choose $\epsilon_{ij}=0.1,~\forall_{\{ij\} \in \mathcal{E}}$ for the transmission line constraints in all settings discussed. For the generation constraints, we choose different $\epsilon_i$ for the different problems, with $\epsilon_{i}=0.001,~\forall_{\{i\} \in \mathcal{G}}$ for the CC-OPF and the linear WCC-OPF and $\epsilon_{i}=10^{-5},~\forall_{\{i\} \in \mathcal{G}}$ for the quadratic WCC-OPF.

The cost of the OPF solutions with the aforementioned $\epsilon_i,~\epsilon_{ij}$ is very similar for the three formulations, as shown in Table \ref{costComparison}. The cost difference is less than 0.1\%, with the chance constraint being slightly cheaper than the other two, and the quadratic weighting function being the most expensive. Since the cost of the different formulations is so similar, it is useful to compare their performance in terms of number and size of violations.

\begin{table}
\caption{Cost for the CC-OPF and the WCC-OPFs}
\label{costComparison}
\centering
\begin{tabular}{|c|c|c|c|}
\hline
                & {CC-OPF}  & {Linear WCC-OPF}  & {Quadratic WCC-OPF} \\
\hline
Total cost [\$] & 16546     &  16547  (+0.01\%) &  16562   (+0.1\%)  \\
\hline
\end{tabular}
\end{table}

To test the OPF solutions, we generate 10 000 samples from the multivariate distribution of the wind power fluctuations and compute the number and size of the resulting violations.
For the generation constraints,
the choice of a very small $\epsilon_i$ ensures that there are very few significant violations of the generation limits.
%the CC-OPF provide a solution with very few violations of the generation limits since $\epsilon_i=0.001$.
%For the WCC-OPF, there are also a negligible amount of violations larger than 0.1 MW.
However, since $\epsilon_i>0$, the linear and quadratic WCC assign a small, but non-zero $\boldsymbol{\alpha}$ to some generators procuding at the limit.
With symmetrically distributed wind power deviations with zero mean, the corresponding generators thus experience a violation probability of 50\%, but with very small violations  $<0.1$ MW.
%lead to very few violations of the generation limit.  a small epsilon lead to very few significant violations. However, since $\epsilon_i>0$, the
%Due to the restrictive choice of $\epsilon_i$, the generation constraints have a very low probability of significant violations. However, since $\epsilon_i$ is non-zero, the weighted chance constraint is able to assign a small $\boldsymbol{\alpha}_i$ even to generators that are already at their upper or lower limit.
%Since the required change in generation output follows the overall wind deviation, which is symmetrically distributed around 0, in half of the cases, these generators will be asked to change their output outside their limits. The resulting violation probability is therefore close to 50\% for those generators. However, since the corresponding $\boldsymbol{\alpha}_i$ are very small, the actual violation is less than 0.1 MW. Therefore, we do not discuss these constraints further, and focus instead on the transmission line constraints.

\begin{figure}[!t]
\includegraphics[width=0.5\textwidth]{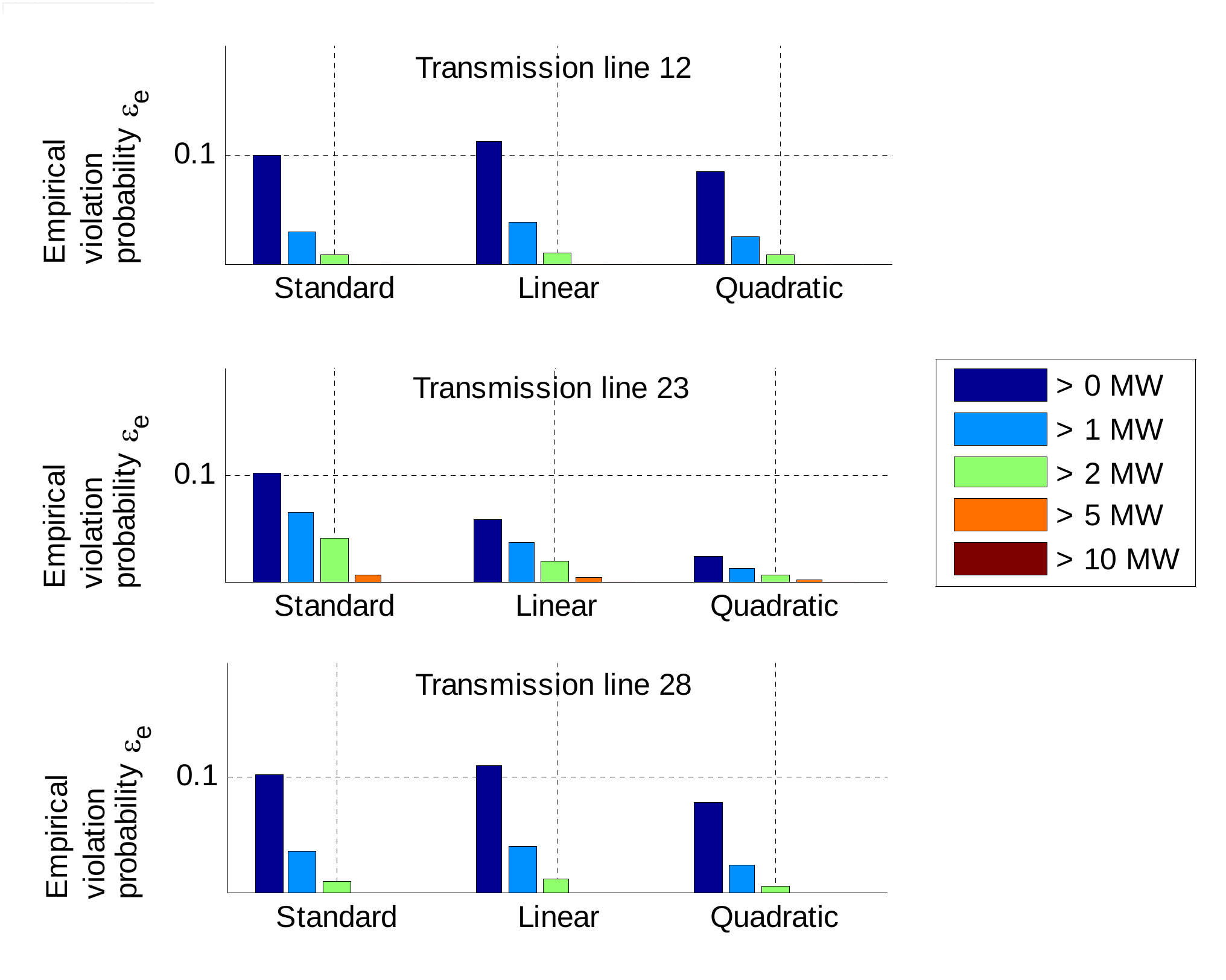}
\centering
\caption{Empirical violation probability ${\epsilon_e}$ for transmission line 12 (upper part), line 23 (middle) and line 28 (bottom). The results for all three formulations (standard, linear and quadratic) are shown in each plot, and the color of the bar indicate the empirical probability of exceeding certain violation thresholds (0, 1, 2, 5 and 10 MW).}
\label{fig_violations}
\end{figure}

For the transmission line constraints, where $\epsilon$ is larger, there are more violations. We only consider the active transmission lines constraints, since this is where the effect of the constraint formulation is visible. The active transmission constraints are also the only ones who experience a significant number of violations.
In Fig. \ref{fig_violations}, the empirical probability ${\epsilon_e}$ of violations larger than 0, 1, 2, 5 and 10 MW is plotted for three transmission lines, line 12 (top), line 23 (middle) and line 28 (bottom).
For the standard CC, the probability of any violation $> 0 MW$ is close to 0.1 for all three lines. The linear WCC lead to a higher violation probability for line 12 and 23, but lower violation probability for line 28. The quadratic WCC is more restrictive, with all violation probabilities below 0.1.
Considering the size of the different violations, we see that line 23 has some violations $> 5 MW$. For this line, it is clear how the standard CC allows for more and larger violations than the linear and quadratic WCC. Particularly the quadratic WCC reduces the number of large overloads to almost zero.

Based on the results, we see that the WCC are less restrictive than the standard CC for small violations, and $\epsilon_e$ exceeds the probability of the corresponding CC. However, both the linear and quadratic WCC are more effective in reducing the probability of large overloads, as seen for line 23. This demonstrates the purpose of enforcing a weighted chance constraint as opposed to a standard chance constraint.

\subsection{Comparison of CC-OPF and WCC-OPF}
We now compare the linear WCC-OPF with an affine policy \eqref{affine} and a piecewise affine policy \eqref{piecewiseaffine} for generation control.
We use the same set-up as above, but choose $\epsilon_{ij}=0.01,~\forall_{\{ij\} \in \mathcal{E}}$ and $\epsilon_{i}=0.01,~\forall_{\{i\} \in \mathcal{G}}$ for both policies. For the piecewise affine policy, we define the threshold for deploying additional reserves $\beta^+,~\beta^-$ as $\Omega^+=-\Omega^-=70~MW$.

\begin{figure}[!t]
\includegraphics[width=0.5\textwidth]{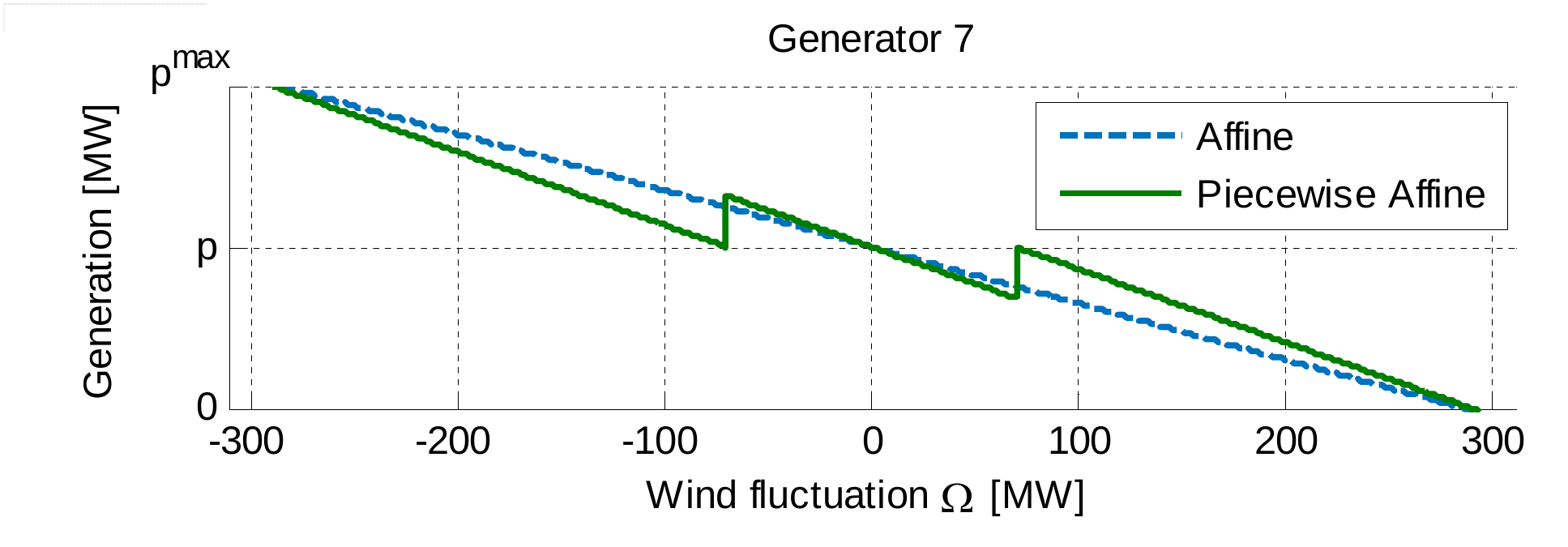}
\centering
\caption{Generation output of generator 7 as a function of the total wind fluctuation $\Omega$ for the affine and the piecewise affine policy.}
\label{fig_generation}
\end{figure}

The cost of the WCC-OPF solutions are shown in Table \ref{comparisonII}. The cost of the CC-OPF with piecewise affine generation control is lower than the policy with affine control.
\begin{table}
\caption{Cost for the affine and piecewise affine policies}
\label{comparisonII}
\centering
\begin{tabular}{|c|c|c|}
\hline
                & {Affine}  & {Piecewise Affine}  \\
\hline
Total cost [\$] & 16546     &  16569  (-0.2\%)   \\
\hline
\end{tabular}
\end{table}
The reason for this cost decrease is due to more flexible use of generation resources. Fig. \ref{fig_generation} shows the generation output of one generator as a function of the total wind fluctuation $\Omega$ for the affine (straight line) and the piecewise affine policy (straight line with jumps).
We observe that both policies lead to similar set-points for the scheduled generation $p$, and that they hit the maximum and minimum generation bounds at the same time. However, the piecewise affine policy allows a steeper slope than the affine policy, due to the use of additional reserves $\beta^+,~\beta^-$ when the total wind fluctuation exceeds the thresholds $\Omega^+,~\Omega^-$.
For the congested transmission lines 12, 23 and 28, a similar effect is observed. For large fluctuations, the additional reserves $\beta^+,~\beta^-$ shifts the power flows to lower values, thus lowering risk.

\section{Conclusions}
In this paper, we introduced the WCC-OPF as an extension to previous CC-OPF formulations.
We defined a new type of chance constraints, the weighted chance constraints, that can handle general control policies while preserving convexity and allowing for efficient computation.
The weighted chance constraints account for the size of constraint violations through a weighting function, which assigns a higher risk to larger violations.
%We prove that the problem remains convex for any convex weighting function, and for very general generation control policies.

The performance of the new WCC-OPF was demonstrated in a case study for the IEEE RTS96 system.
It was shown that the WCC-OPF allows for a larger number of small violations, but reduces the number of severe overloads compared to the CC-OPF.
Further, %the original affine control policy for generation control was compared with a piecewise affine policy which models the deployment of slow reserve capacities.
it was shown that a more flexible, piecewise affine policy for generation control allows for lower cost compared with the original affine policy, while maintaining the same risk level.

%To compare the CC-OPF with the WCC-OPF, all chance constraints were assigned the same type of weighting function. It would however be possible to use different weighting functions for generation and transmission constraints, since the consequences of the different types of violations are dissimilar.
In our computations, both transmission and generation constraints were assigned the same type of weighting function. It would however be possible to use different weighting functions for generation and transmission, since the consequences of a violation is different for the two.
This type of analysis, along with an extension of the standard chance constraints to piecewise affine policies, is part of future work.

%{\color{red}
Convexity of the WCC-OPF formulation discovered in this manuscript suggests a path forward for scaling up the approach to larger systems through deployment and development, in the spirit of \cite{12BCH}, of a variety of computational gadgets of the optimization theory such as cutting plane, piece-wise linear approximations of convex constraints and others. Notice also that our method allows generalizations to account for non-Gaussian effects in fluctuations of renewables,  while preserving convexity of the formulation. Future work is needed to explore possible synergy of this observation with complementary work on improving statistical modeling of renewables \cite{Yury} and further exploration of advanced sampling techniques \cite{kostas}, \cite{george}.
The possibility to formulate general control policies can also be exploited to flexibly model other power system controls, such as FACTS or HVDC, extending the current affine models applied in \cite{UMBRELLA}.
Finally, we plan combining collection of ideas and techniques proposed in the manuscript with recent ideas based on energy function \cite{Misha-IREP13,Dj-Energy} and monotone operator \cite{Dj-Monotone} approaches to attack more challenging CC- and WCC versions of the nonlinear AC-OPF formulations, thus accounting for risks beyond the DC-approximation explored so far.
%}

%%%%%%%%%%%%%%%%%%%%%%%%%%%%%%%%%%%%%%%%%%%%%%%%%%%%%%%%%%%%%%%%%%%%%%%%%%%%%%%%
\section{Acknowledgements}

%{\color{red}
The authors thank S. Backhaus and D. Bienstock for multiple discussions and advice.
The work at LANL was carried out under the auspices of the National Nuclear Security Administration of the U.S. Department of Energy at Los Alamos National Laboratory under Contract No. DE-AC52-06NA25396 and it was partially supported by DTRA Basic Research Project $\#10027-13399$. The authors also acknowledge partial support of the Advanced Grid Modeling Program in the US Department of Energy Office of Electricity.
%}
Line Roald receives funding from the project "Innovative tools for future coordinated and stable operation of the pan-European electricity
transmission system" (UMBRELLA), supported under the 7th Framework Programme of the European Union, grant agreement 282775.

%The authors gratefully acknowledge the contribution of National Research Organization and reviewers' comments.

%%%%%%%%%%%%%%%%%%%%%%%%%%%%%%%%%%%%%%%%%%%%%%%%%%%%%%%%%%%%%%%%%%%%%%%%%%%%%%%%

\bibliographystyle{IEEEtran}
\bibliography{20150223_bib_LR}
%References are important to the reader; therefore, each citation must be complete and correct. If at all possible, references should be commonly available publications.
%
%\begin{thebibliography}{99}
%
%\bibitem{c1}
%J.G.F. Francis, The QR Transformation I, {\it Comput. J.}, vol. 4, 1961, pp 265-271.
%
%\bibitem{c2}
%H. Kwakernaak and R. Sivan, {\it Modern Signals and Systems}, Prentice Hall, Englewood Cliffs, NJ; 1991.
%
%\bibitem{c3}
%D. Boley and R. Maier, "A Parallel QR Algorithm for the Non-Symmetric Eigenvalue Algorithm", {\it in Third SIAM Conference on Applied Linear Algebra}, Madison, WI, 1988, pp. A20.
%
%\end{thebibliography}

\end{document}